\documentclass[
a4paper,
envcountsame,
linenumbering
]{llncs}

\usepackage[utf8]{inputenc}
\usepackage[textsize=scriptsize,obeyFinal]{todonotes}

\usepackage[T1]{fontenc}
\usepackage{color}
\usepackage{xcolor}

\usepackage{lineno}
\linenumbers

\usepackage{amsmath, amsfonts, amssymb, dsfont}

\usepackage{enumerate}

\usepackage{wrapfig}
\usepackage{caption} 
\usepackage{float}

\newtheorem{observation}[theorem]{Observation}

\def\ta{\mathtt{a}}

\renewcommand{\epsilon}{\varepsilon}
\renewcommand{\phi}{\varphi}

\def\N{\mathbb{N}} 
\def\Z{\mathbb{Z}} 

\DeclareMathOperator{\Fact}{Fact}
\DeclareMathOperator{\Pref}{Pref}
\DeclareMathOperator{\pref}{pref}
\DeclareMathOperator{\Suff}{Suff}

\DeclareMathOperator{\Parents}{Parents}
\DeclareMathOperator{\Children}{Children}
\DeclareMathOperator{\Factone}{Fact1}

\usepackage{bbm}

\def\ta{\mathtt{a}}

\def\hamd{\operatorname{d}_{\operatorname{H}}}
\DeclareMathOperator{\PN}{PN}
\def\parpn{\Parents_{\PN}}

\def\N{\mathbb{N}} 
\def\Z{\mathbb{Z}} 

\def\pref{\operatorname{pref}}
\def\suf{\operatorname{suf}}

\renewcommand{\epsilon}{\varepsilon}
\renewcommand{\phi}{\varphi}

\DeclareMathOperator{\flip}{flip}

\usetikzlibrary{positioning, shapes,patterns,shadows.blur, arrows,decorations,arrows,automata,shadows,patterns,chains,graphs,calc,intersections,matrix,fit,shapes,chains,decorations.pathreplacing}
\usepackage{paralist}
\usepackage{algorithm2e}

\usepackage{booktabs}
\usepackage{float}
\usepackage{hyperref}

\usepackage{cleveref}

\definecolor{cau}{RGB}{156,10,125} 

\newif\ifpaper
\paperfalse 

\nolinenumbers

\title{Word Chain Generators for Prefix Normal Words}

\titlerunning{Word Chain Generators}

\author{Duncan Adamson\inst{1}\and Moritz Dudey\inst{2}\and
Pamela Fleischmann\inst{2} \and
Annika Huch\inst{2}}

\institute{University of St Andrews, UK\\
	\email{duncan.adamson@st-andrews.ac.uk}\and Kiel University, Germany, \email{stu227171@mail.uni-kiel.de, $\{$fpa,ahu$\}$@informatik.uni-kiel.de}}

\authorrunning{D. Adamson, \and M. Dudey, \and P. Fleischmann, \and A. Huch}

\begin{document}

\maketitle

\begin{abstract}
In 2011, Fici and Lipták introduced prefix normal words. A binary word is prefix normal if it has no factor (substring) that contains more occurrences of the letter 1 than the prefix of the same length. Among the open problems regarding this topic are the enumeration of prefix normal words and efficient testing methods. We show a range of characteristics of prefix normal words. These include properties of factors that are responsible for a word not being prefix normal. With word chains and generators, we introduce new ways of relating words of the same length to each other.
\end{abstract}

\section{Introduction}
A binary word is \emph{prefix normal} if it contains no factor (substring) with more occurrences of the letter 1 than the prefix of the same length \cite{DBLP:conf/dlt/FiciL11}. For example, 1101 is prefix normal, as every factor contains no more 1s than the prefix of the same length, while 100101 is not prefix normal since its factor 101 contains more 1s than the prefix of the same length, 100.
Research on prefix normal words emerged from work on \emph{binary jumbled pattern matching} (BJPM) \cite{DBLP:journals/ijfcs/BurcsiCFL12,DBLP:journals/algorithmica/GagieHLW15}. This is a combinatorial problem asking, given a word over the binary alphabet $\{0,1\}$ and two natural numbers $x$ and $y$, whether the word contains a factor with $x$ 0s and $y$ 1s. Note that the order in which these letters appear is considered irrelevant (hence the name \emph{jumbled} pattern matching). The description of a word by the number of occurrences of each letter in the underlying alphabet is given by its \emph{Parikh vector} \cite{DBLP:journals/jacm/Parikh66,DBLP:journals/eatcs/Salomaa03}. If we look at all factors of length $x + y$ of a word (sliding window), we can answer a BJPM query in linear time. The \emph{indexed} jumbled pattern matching problem (e.g.,~\cite{DBLP:conf/spire/LeeLZ12}) is concerned with answering multiple such queries on the same word. One way to construct an index is by computing the maximum amount of 1s and 0s in any factor of any length. With that, we can create a new word where each prefix contains exactly the maximum number of 1s in any factor of the same length of the original word. In \cite{DBLP:conf/dlt/FiciL11}, this new word was introduced as the \emph{prefix normal form} of a word with respect to 1. For example, for 00101 we will find that the word 10100 fulfils this property. 
Thus, with the prefix normal form of a word this maximum can be looked up in linear time, or for multiple searches a look-up table can be computed in linear time such that the maximum can be obtained in $O(1)$.
As shown in \cite{DBLP:conf/stringology/CicaleseFL09}, for any number between the maximum and the minimum of 0s in a word, the word always contains a factor with that number of occurrences of 0s. That explains why prefix normal forms can be used to create an index for BJPM queries \cite{DBLP:conf/dlt/FiciL11}. Calculating the prefix normal form of a word in such a manner requires a quadratic amount of time, and indeed provides one of the main motivations for studying prefix normal words: A more efficient way to generate prefix normal forms would lead to a faster solution for the indexed BJPM problem. As presented in \cite{DBLP:journals/ijfcs/BurcsiCFL12}, all BJPM queries for a word can be answered in $O(1)$ if its prefix normal form is known.

%
In order to obtain the prefix normal form of a word, Fici and Lipták introduced in \cite{DBLP:conf/dlt/FiciL11} the notion of \emph{prefix normal equivalence}. A word is \emph{prefix normal equivalent} to another if both have the same prefix normal form. The proof that this is indeed an equivalence relation was presented by Burcsi et al. in \cite{DBLP:journals/tcs/BurcsiFLRS17}. They also showed that each equivalence class has a unique prefix normal representative. Fici and Lipták's work \cite{DBLP:conf/dlt/FiciL11} includes a characterization of prefix normality using two functions: For a binary word $w$, $f_w$ maps a natural number $x$ to the maximum number of 1s in any factor of length $x$ in $w$. The other function, $p_w$, maps a natural number $y$ to the amount of 1s in the prefix of $w$ of length $y$. A word $w$ is prefix normal if $f_w$ equals $p_w$. 
Further, Fici and Lipták proved that the language of prefix normal words is not context-free, and explored the relationship between Lyndon and prefix normal words, showing that every prefix normal word is a pre-necklace.
Among the open problems stated are the unknown size and number of prefix normal equivalence classes.
It was shown in \cite{DBLP:conf/cpm/BurcsiFLRS14} that prefix normal words form a \emph{bubble language}. This means that the first occurrence of $01$ in a prefix normal word can be swapped to $10$, resulting in another prefix normal word. This property has been used in new algorithms for testing \cite{DBLP:conf/fun/BurcsiFLRS14} and an algorithm generating all prefix normal words of a fixed length \cite{BURCSI202086}. In \cite{DBLP:conf/lata/FleischmannKNP20}, prefix normal palindromes and \emph{collapsing} prefix normal words were examined. Two words $v$ and $w$ of the same length \emph{collapse} if $1v$ and $1w$ are prefix normal equivalent. \emph{Infinite} prefix normal words were studied in \cite{DBLP:journals/tcs/CicaleseLR18,DBLP:journals/tcs/CicaleseLR21}. In \cite{DBLP:journals/tcs/BurcsiFLRS17}, it was conjectured that the number of prefix normal words of length $n$ is bounded by $2^{n-\Theta((\log n)^2)}$. This has since been proven in \cite{BALISTER201975}. A closed form formula for the number of prefix normal words or a generating function has not been found yet.

\begin{figure}
	\includegraphics[width=0.5\textwidth]{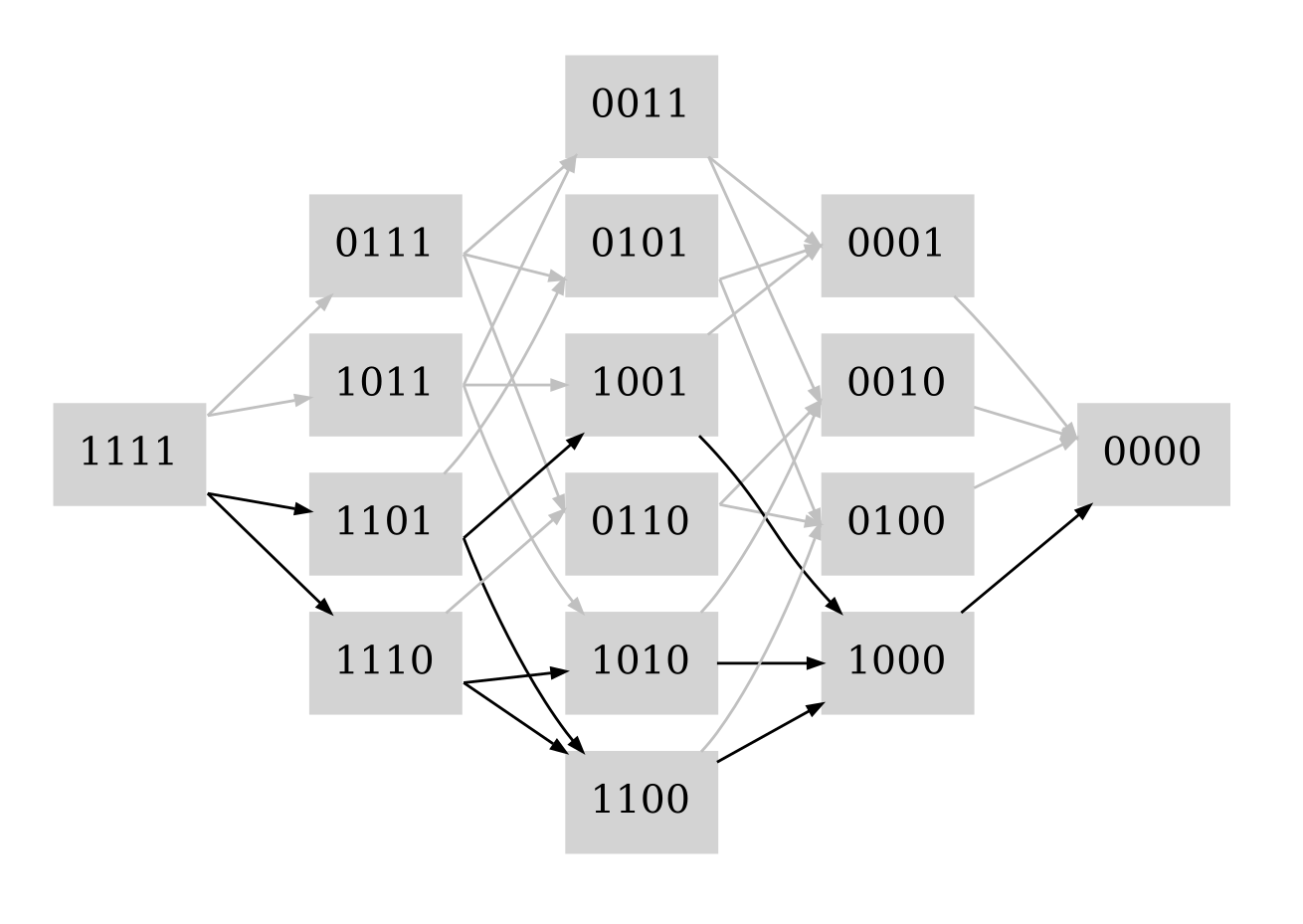}
	\centering
	\caption{Graph of words of length 4 with paths corresponding to prefix normal words in black.}
	\label{fig:graph}
\end{figure}

\textbf{Our contribution.} 
In our work, we start by presenting observations on (minimal) factors of prefix normal words. We point out several characteristics of factors that are responsible for a word not being prefix normal. 
Our approach to the enumeration of prefix normal words, comes with the new concept of \emph{prefix normal word chains} and \emph{prefix normal word chain generators}. Here, we observe graphs of binary words of length $n$ where every node is labelled with some $w \in \Sigma^n$ and two nodes are adjacent if their two words differ in only one letter (\Cref{fig:graph}). We investigate the paths starting with $1^n$ and reaching $0^n$ within those graphs that only contain prefix normal words (marked non-transparent in \Cref{fig:graph}). Since every path corresponds to a permutation of $\{1, \ldots, n\}$ that indicated at which position a letter is changed from $1$ to $0$, we introduce those special permutations as \emph{prefix normal word chain generators} whereas we refer to the respective sequence of prefix normal words along this path as \emph{prefix normal word chain}.
After some basic observations, we examine under which circumstances two numbers in a prefix normal generator can be swapped to construct a new such generator. Then, we connect prefix normal generators with the iterative construction of prefix normal words introduced in 
\cite{DBLP:conf/fun/BurcsiFLRS14,BURCSI202086,DBLP:conf/lata/FleischmannKNP20}.

\textbf{Structure of the work.} \Cref{prelims} is about basic definitions. Results on more efficient testing on prefix normality can be found in \Cref{minFactors}. \Cref{wordchains} is concerned with the enumeration of prefix normal words. There, we introduce \emph{word chains} and \emph{word chain generators}. 

\section{Preliminaries}\label{prelims}
\label{sec:prelims}
	Let $\N$ be the set of natural numbers starting with 1. Let $\N_0 = \N \cup \{0\}$. For $n \in \mathbb{N}$ let $[n] = \{i \in \mathbb{N} \: | \:  i \leq n\}$ and $[n]_0 = [n] \cup \{0\}$. For a sequence $s$, let $s[k]$ be the $k^{\text{th}}$ element of that sequence. 

An \emph{alphabet} is a finite set of elements called \emph{letters}. We will often use the \emph{binary alphabet} $\Sigma = \{0,1\}$.
A \emph{word} is a sequence of letters. Let $\Sigma^*$ denote the set of all finite words consisting of letters in $\Sigma$. For $w \in \Sigma^*$, $|w|$ is the \emph{number of letters} in $w$ and $w[i]$ is the $i^{\text{th}}$ letter in $w$ for $i \in [|w|]$. Let $\varepsilon$ denote the \emph{empty word}, meaning $|\varepsilon| = 0$. 
Let $x,y,z \in \Sigma^*$ with $w = xyz$, then  $y \in \Sigma^*$ is called a {\em factor}. Let $\text{Fact}(w)$ denote the set of factors of a word $w$, and for $k \in \mathbb{N}$, let $\text{Fact}_k (w)$ be the set of factors of length $k$ of $w$. 
Further, $x$ is a \emph{prefix} of $w$, which is denoted as $x \leq_p w$, and $z$ is a \emph{suffix} of $w$, denoted as $z \leq_s w$.  Let for any $k \in [|w|]$, $\pref_k(w)$ denote the prefix of length $k$ of a word $w$, and $\text{suf}_k(w)$  the suffix of length $k$ of a word $w$. 
Let $w[i:k]$ denote the factor $w[i] \cdots w[k]$, for $i, k \in [|w|]$ with $i \leq k$. If $i > k$ or, for some $i,k \in \Z$, $i \not\in [|w|]$ or $k \not\in [|w|]$, then $w[i:k] = \varepsilon$. 
For $v \in \Sigma^*$ and $w \in \Sigma^*$, we define $|w|_v = |\{(i,j) \in [|w|]^2 \: | \: w[i:j] = v\}|$, i.e. $|w|_v$ denotes the number of times $v$ occurs as a factor in $w$. Two factors $x$ and $y$ of a word $w \in \Sigma^*$ \emph{overlap} with each other if they share at least one index of $w$, i.e. if there exist $i,j,k,l \in \N$ with $x = w[i:j], y = w[k:l]$ and $i \leq j, k \leq l$, then $x$ and $y$ overlap if $\{i,i+1,...,j\} \cap \{k, k+1, ..., l\} \neq \emptyset$. We say that for $w \in \Sigma^*$ and a factor $v = w[k:l]$ for $k,l \in [|w|]$, $v$ \emph{includes} an index $a \in [|w|]$ of $w$ if $k \leq a \leq l$. 
A word $w$ is called a \emph{palindrome} if $w = w^R$ where $w^R$ denotes the word in reverse. 
For a word $w \in \Sigma^*$ and $k \in \N_0$, let $w^0 = \varepsilon$ and $(w)^{k+1} = ww^k$. 
A \emph{block} is a maximal unary factor of a word, i.e. $u\in\Fact(w)$ is a block if there exists $\ta\in\Sigma$ such that $u=\ta^{|u|}$ and either $u\leq_p w$ and $w[|u|+1]\neq \ta$ or $u\leq_s w$ and $w[|w|-|u|]\neq\ta$ or there exist
$x,y\in\Sigma^+$ with $w=xuy$ and $x[|x|]\neq\ta\neq y[1]$. 
Blocks of 0s and 1s are called 0- and 1-runs, resp., in \cite{BURCSI202086}.

A \emph{permutation} of a set $S$ is a bijection $\sigma$ from $S$ to $S$. We call $|S|$ the \emph{length of the permutation}. 
As an example, for the set $S = \{1,2,3,4\}$, the bijection $\sigma : S \to S, x \mapsto (x+1) \text{ mod } 4 \:$ is a permutation, and the bijection $\gamma : S \to S$, $\gamma (1) = 4, \gamma(2) = 2, \gamma(3) = 1, \gamma(4) = 3 \:$ is also a permutation. We will write these as words $\sigma = 2341$ and $\gamma = 4213$.

The \emph{Hamming distance} between two words of equal length is the number of indices at which the two words differ: for $u,v \in \Sigma^*$, with $|u| = |v|$, $\hamd(u,v) = |\{i \: | \: i \in [|u|], x[i] \neq y[i]\}|$ denotes the Hamming distance between $u$ and $v$. 

\medskip

After these basic definitions, we introduce the main object of interest - the prefix normal words.

\begin{definition}
For $w \in \Sigma^*$, define the \emph{prefix-ones function} and the \emph{maximum-ones function} by
$p_w : [|w|] \to [|w|], i \mapsto |\pref_i(w)|_1,
f_w : [|w|] \to [|w|], i \mapsto \max_{u \in \Fact_i(w)} |u|_1$.
\end{definition}

\begin{definition} 
A word $w \in \{0,1\}^*$ is called \emph{prefix normal} if and only if $p_w = f_w$. If $w$ is not prefix normal there exists an $i\in[|w|]$ and $x\in\Fact_i(w)$ with $|x|_1>\pref_i(w)$. We say that {\em $x$ is responsible for $w$ not being prefix normal}.
\end{definition}

For instance, the word $w = 111001011011$ is not prefix normal, because its factor $w[8:12] = 11011$ contains more 1s than the prefix of length five, namely $11100$. The prefix normal words up to length three are: $0$, $1$, $00$, $10$, $11$, $000$, $100$, $101$, $110$, $111$.


%


For the definition of a word chain, we need a function that takes a word and flips a single $1$ to a $0$. For example,  $\flip_1(001) = 001$ while $\flip_2(11010) = 10010$.  A similar operation that allows for flipping $0$s to $1$s was defined in \cite{DBLP:journals/tcs/CicaleseLR18}. 

\begin{definition}
    Let $i \in \mathbb{N}$. Define $\operatorname{flip}_i : \Sigma ^* \to \Sigma^*, w \mapsto w[1:i-1] \cdot 0 \cdot w[i+1:|w|]$, if $1 \leq i \leq |w|$ and $w \mapsto w$, otherwise.
\end{definition}


\begin{definition}
    For two words $p,c \in \Sigma^*$, we call $p$ a \emph{parent} of $c$ and $c$ a \emph{child} of $p$ if and only if 
       $|p| = |c|$,
       $|p|_1 - 1 = |c|_1$ and
       $\hamd(p,c) = 1$.
    Let $\Parents(w)$ denote the \emph{set of all parents} of a word $w \in \Sigma^*$ and let $\Children(w)$ denote the \emph{set of all of its children}. Let $\parpn(w) = \{x \in\Sigma^{\ast} |\,  x \in \Parents(w) \land x \mbox{ prefix is normal}\}$.
\end{definition}

\begin{remark}
Another way of describing this relation of parents and children is by using the $\flip$ function defined above. Let $w \in \Sigma^*$. Then for all $i \in [|w|]$, if $w[i]=1$ then $w$ is a parent of $\flip_i(w)$ and $\flip_i(w)$ is a child of $w$.
\end{remark}

Our last definition introduces the prefix normal word chains and their generators. The idea is to start with the prefix normal word $1^k$ for some $k\in\N$ and flipping successively $1$s to $0$s in such an order that every word in the chain (until we reach $0^k$) is prefix normal. This sequence of prefix normal words is given by a permutation of the word's indices.

\begin{definition}
    Let $\sigma$ be a permutation on $[m]$ for some $m \in \N$. Let the sequence of words emerging from applying the flip function at the indices given by $\sigma$ be 
    $$
    c_{\sigma} = \left( \: \flip_{\sigma(n-1)} (1^m) \circ\cdots\circ  \flip_{\sigma(1)} (1^m)\right) _{n \in [m+1]}.
    $$
    Such a sequence $c_{\sigma}$ is a \emph{word chain} and $\sigma$ its \emph{word chain generator}.
    We call $c_\sigma$ a \emph{prefix normal word chain} and $\sigma$ the corresponding \emph{prefix normal word chain generator} (or for short {\em prefix normal generator}) iff for all $n \in \mathbb{N}$, $c_{\sigma}[n]$ is prefix normal. Here, we denote by $c_{\sigma}[n]$ the projection onto the $n^{\text{th}}$ component of $c_{\sigma}$.
\end{definition}

Note that word chains can also be defined recursively, i.e., $c_\sigma[1] = 1^m$ and $c_\sigma[n+1] = \flip_{\sigma(n)}(c_\sigma[n])$ for a permutation $\sigma$ on $[m]$ for some $m,n \in \N$, $n\leq m$.

The sequence $\sigma = 231$ is a prefix normal word chain generator, because it creates the prefix normal word chain $c_{\sigma} = (111, 101, 100, 000)$. The permutation $\gamma = 213$ generates the sequence $s_{\gamma} = (111, 101, 001, 000)$ and is not a prefix normal generator as $s_\gamma[3] = 001$ is not a prefix normal word.



\subsection{Observations on (Minimal) Factors of Prefix Normal Words}\label{minFactors}
\label{sec:minFactors}

In this section we present observations about prefix normal words. The first idea is to characterize the minimal factor that contains more 1s than the prefix of the same length.

\begin{observation}{}
    \label{lem:beginAndEndIn1}
    The factor minimal in length that is responsible for a word not being prefix normal always begins and ends with a 1.
\end{observation}
\ifpaper 
\else 
    \begin{proof}
        Let $w$ be a binary word and let $u$ be the shortest factor of $w$ with $|u|_1 > |\pref_{|u|}(w)|_1$. Now assume that $u$ begins (ends) with a 0. Let $v$ be the suffix (prefix) of $u$ without that 0, so $v = \suf_{|u|-1}(u)$ (or $v=\pref_{|u|-1}(u)$ resp.). Because the first (last) letter of $u$ is a 0, and all other letters are the same, we know that $|u|_1 = |v|_1$. Since a shorter word cannot contain more 1s than a longer one, we also know that $|\pref_{|u|}(w)|_1 \geq |\pref_{|u|-1}(w)|_1$. That gives us $$|v|_1 = |u|_1 > |\pref_{|u|}(w)|_1 \geq |\pref_{|u|-1}(w)|_1 = |\pref_{|v|}(w)|_1.$$ This is a contradiction to $u$ being the shortest factor of $w$ responsible for $w$ not being prefix normal. So this shortest factor cannot begin (end) with a 0.\qed
    \end{proof}

\fi

\begin{observation}{}
    \label{lem:minFacOneMoreOne}
    The factor minimal in length responsible for a word not being prefix normal always contains just one more 1 than the prefix of the same length.
\end{observation}
\ifpaper 
\else 
    \begin{proof}
        Let $w \in \Sigma^*$, and let $x \in \Fact(w)$ and $y = \pref_{|x|}(w)$ such that $x$ is the minimal factor with $|x|_1 > |y|_1$. Now assume $|x|_1 > |y|_1 + 1$. This, however, means that the factor $x' = \pref_{|x|-1}(x)$ also contains more ones than $y' = \pref_{|x'|}(y)$, because $|x'|_1 \geq |x|_1 - 1 > |y|_1 \geq |y'|_1$. This is a contradiction to $x$ being minimal. So $x$ can only contain exactly one extra 1 compared to $y$.\qed
    \end{proof}

\fi

By Observation~\ref{lem:beginAndEndIn1}, we only have to check all factors in a word that both begin and end with a 1, which reduces the complexity compared to checking every factor. For $w \in \Sigma^*$, the number of factors that $w$ contains is $\frac{|w| \cdot (|w|+1)}{2}$. Note that $|\Fact(w)| \leq \frac{|w| \cdot (|w|+1)}{2}$, because factors can occur more than once. 
With the following result, we give the number of factors within a word starting and ending with a $1$, including the trivial factor $1$.

\begin{remark}{}
    Let $w\in\Sigma^{\ast}$. Let $\Factone(w)$ denote the set of factors of $w$ both starting and ending in $1$. Let 
    $f^{1...1}(w) = \sum_{u \in \Factone(w)} |w|_u$
    denote the number of factors in $w$ starting and ending with $1$. Then 
    $$f^{1...1}(w) = \frac{|w|_1 (|w|_1 - 1)}{2}.$$
\end{remark}


Before we determine which characteristics at least one factor of a non prefix normal word must have, we need the following observation about prefix normal words in general.

\begin{observation}{}
    \label{lem:FacInPrefixOfOnes}
    A factor that includes any part of the first block of 1s of a word cannot be responsible for the word not being prefix normal. Let $w \in \Sigma^*, k \in \N, u \in \Sigma^*$ such that $w = 1^ku$ and a factor $x = w[i:j], i,j \in [|w|], i \leq j$. If $i \leq k$, then $|x|_1 \leq |\pref_{|x|}(w)|_1$.
\end{observation}
\ifpaper 
\else    
    \begin{proof}
        We know that $|w[1:i-1]| = |w[1:i-1|_1$, i.e., $w$ consists only of 1s up to the index $i$. So $w[i-1] = 1$. Now, compare $x$ to the factor of the same length, shifted one index to the left. This factor's first letter, $w[i-1]$, is a 1 that $x$ does not contain. On the other side, the last letter of $x$, namely $w[j]$ is either $1$ or $0$ which this factor does not contain. Thus, $|w[i-1:j-1]|_1 \geq |w[i:j]|_1 = |x|_1$. And since there are only 1s to the left of $w[i]$, this applies to every factor on the left of $x$. So for any number $\ell \in [i-1]$, $|w[i-\ell-1:j-\ell-1]|_1 \geq |w[i-\ell:j-\ell]|_1$ holds. We have       \begin{align*}
        & \: |\pref_{|x|}(w)|_1 \\
        = & \: |w[1:|x|]|_1 \\
        = & \: |w[1:j-i+1]|_1 \\
        \geq & \: |w[2:j-i+2]|_1 \\
        \geq & \: ... \\
        \geq & \: |w[i:j-i+i]|_1 \\
        = & \: |w[i:j]|_1 \\
        = & \: |x|_1,
        \end{align*}
        and therefore, $x$ cannot be a factor that is responsible for $w$ not being prefix normal.\qed
    \end{proof}

\fi

Informally, if a factor \emph{is} responsible for a word not being prefix normal, meaning it \emph{does} contain more 1s than the prefix of the same length, then it does \emph{not} begin inside the first block of 1s in the word. With this observation, we can now show the following proposition that extends Observation~\ref{lem:beginAndEndIn1} that does not refer to blocks of $1$s.

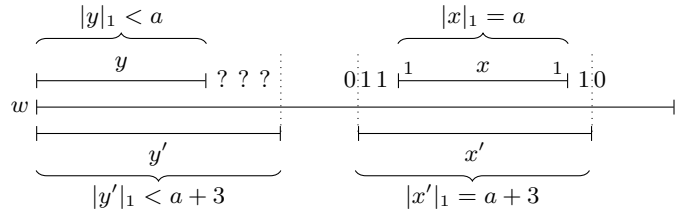
\begin{figure}[h!]
    \centering
        \begin{tikzpicture}[scale=0.7, decoration={brace, amplitude=4pt}]
            \node at (1.6, 1.7) {$|y|_1 < a$};
            \node at (8.4, 1.7) {$|x|_1 = a$};
            \draw[decorate] (0,1.2) -- (3.2,1.2);
            \draw[decorate] (6.8,1.2) -- (10,1.2);
            \draw[|-|] (0,0.5) -- (3.2,0.5) node[pos=0.5, sloped, above] {$y$};
            \node at (7,0.7) {\scriptsize $1$\strut};
            \node at (9.8,0.7) {\scriptsize $1$\strut};
            \draw[|-|] (6.8,0.5) -- (10,0.5) node[pos=0.5, sloped, above] {$x$};

            \node at (-0.3, 0) {$w$};
            \draw[|-|] (0,0) -- (12,0);
            
            \node at (3.5,0.5) {$?$\strut};
            \node at (3.9,0.5) {$?$\strut};
            \node at (4.3,0.5) {$?$\strut};
            \node at (5.9,0.5) {$0$\strut};
            \node at (6.2,0.5) {$1$\strut};
            \node at (6.5,0.5) {$1$\strut};
            \node at (10.3,0.5) {$1$\strut};
            \node at (10.6,0.5) {$0$\strut};

            \draw[dotted] (4.6,1) -- (4.6, -0.5);
            \draw[dotted] (6.05,1) -- (6.05, -0.5);
            \draw[dotted] (10.45,1) -- (10.45, -0.5);
            
            \draw[|-|] (0,-0.5) -- (4.6,-0.5) node[pos=0.5, sloped, below] {$y'$};
            \draw[|-|] (6.05,-0.5) -- (10.45,-0.5) node[pos=0.5, sloped, below] {$x'$};

            \draw[decorate, decoration={mirror}] (0,-1.2) -- (4.6,-1.2);
            \draw[decorate, decoration={mirror}] (6,-1.2) -- (10.4,-1.2);            
            \node at (2.3, -1.7) {$|y'|_1 < a+3$};
            \node at (8.2, -1.7) {$|x'|_1 = a+3$};
        \end{tikzpicture}
    \caption{Illustration of Proposition \ref{prop:FactorStartsEndsWithBlock}}
    \small
    \label{fig:StartAndEndInBlockOf1s}
\end{figure}

\begin{proposition}{}
    \label{prop:FactorStartsEndsWithBlock}
    In any word that is not prefix normal, there is a factor that starts and ends with a block of 1s and contains more 1s than the prefix of the same length (see Figure \ref{fig:StartAndEndInBlockOf1s}). 
\end{proposition}
\ifpaper 
\else    
    \begin{proof}
        Let $w \in \Sigma^*$ and $x \in \Fact(w)\backslash\Pref(w)$,  $y = \pref_{|x|}(w)$ its corresponding prefix, and let $x$ be minimal in length such that $|x|_1 > |y|_1$. If $x$ already begins and ends with a block of 1s, we are done. Somewhere to the left of $x$ in $w$, there has to be a 0 (see Lemma \ref{lem:FacInPrefixOfOnes}). If we extend $x$ to the left and right in a way that only additional 1s are included, the resulting factor $x'$ has been increased by $k$ 1s, where $k = |x'| - |x|$. Let $y' = \pref_{|x'|}(w)$. Now we have $|x'|_1 = |x|_1 + k > |y|_1 + k \geq |y'|_1$.\qed
    \end{proof}

\fi

Note that in Figure \ref{fig:StartAndEndInBlockOf1s}, $x$ is the minimal factor containing 
more 1s than the prefix of the same length. 

\section{Word Chains and Generators}\label{wordchains}
\label{sec:wordchains}
In this section, we provide insights to the problem of enumerating prefix normal words. We approach this by restricting our field of view to just the prefix normal words of fixed length. This leads to \emph{prefix normal word chains} and \emph{prefix normal word chain generators}. First, we give some results on these chains before we connect them to the known field of {\em extension-critical} prefix normal words. 
Note that for any permutation $\sigma$, $c_{\sigma}[1] = 1^m$, because the flip function is applied 0 times.  The following remarks both follow immediately from the definition.

\begin{figure}
    \includegraphics[width=0.5\textwidth]{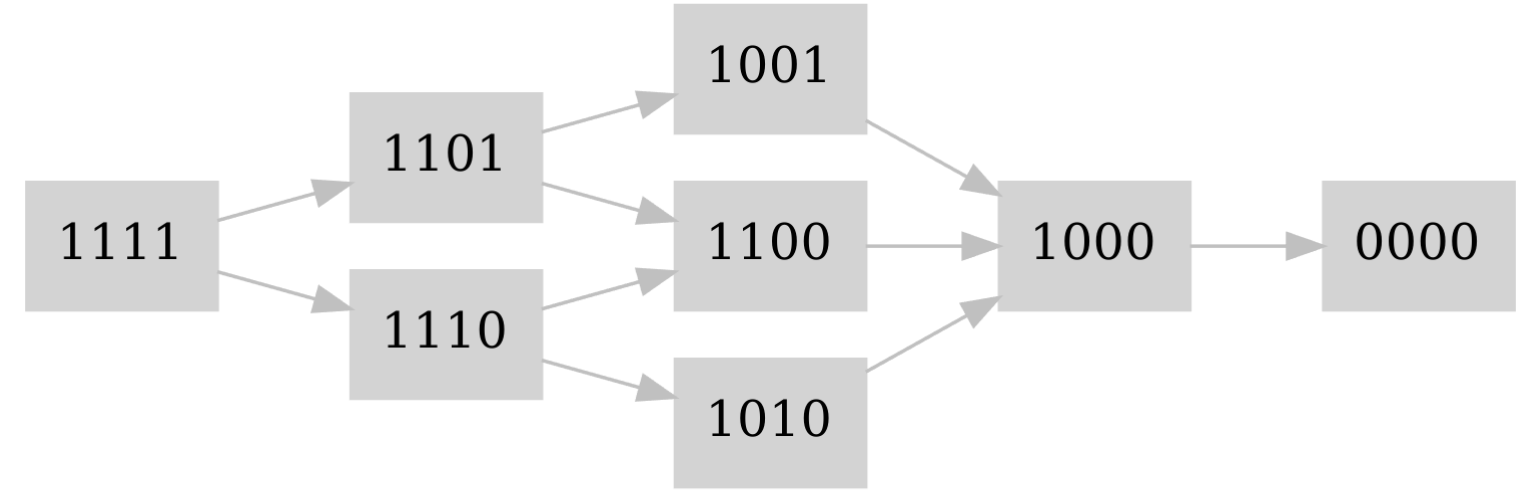}
    \centering
    \caption{Prefix normal words of length 4}
    \label{fig:pn3}
\end{figure}

\begin{remark}{}
    Let $n \in \N$, $\sigma$ a permutation of length $n$, $c_\sigma$ the corresponding word chain and $i \in [n]$. Then, there exist $u,v \in \Sigma^*$ with $|u| = \sigma(i)-1$ and $|v| = n - \sigma(i)$ such that $c_\sigma[i] = u1v$ and $c_\sigma[i+1] = u0v$.  
    \label{rem:zeroAtIndex}
\end{remark}

\begin{remark}{}
    \label{rem:numZerosInChain}
    Let $n \in \N$, $\sigma$ a permutation of length $n$ and $c_\sigma$ the corresponding word chain. For any $i \leq n$, we have $|c_\sigma[i+1]|_0 = i$. 
\end{remark}

Every prefix normal word except those in $\{1\}^*$ has at least one prefix normal parent, since the first occurrence of $0$ can always be replaced by a $1$. With this in mind, the following remark is clear.

\begin{remark}{}
    \label{rem:allpnwInChain}
    Every prefix normal word is part of a prefix normal word chain: every word is part of a word chain and thus, by definition every prefix normal word is part of a prefix normal word chain.
\end{remark}

Every word in a prefix normal word chain except the last one begins with $1$. The index $1$ can only be flipped if all other letters are $0$s.

\begin{lemma}
    Let $n \in \N$ and let $\sigma$ be a permutation over $[n]$. If $\sigma$ is a prefix normal word chain generator, then $\sigma(n) = 1$. 
    \label{lem:1ComesLastInGen}
\end{lemma}
\ifpaper 
\else    
    \begin{proof}
        By contraposition. Let $n \in \N$ and let $\sigma$ be a permutation over $[n]$. Let $s_\sigma$ be the sequence of words created by $\sigma$. Assume that $\sigma(n) \neq 1$. Then therre exists $k < n$  such that $\sigma(k) = 1$. We know from Remark \ref{rem:zeroAtIndex} that $s_\sigma[k+1] = 0w$ for some $w \in \Sigma^*$ with $|w| = n-1$. And Remark \ref{rem:numZerosInChain} gives us 
        \[
        |w|_0 = |0w|_0 - 1 = |s_\sigma[k+1]|_0 - 1 = k - 1. 
        \]
        So $|w|_1 = n - 1 - |w|_0 = n - 1 - (k - 1) > 0$ because $k < n$. Since at least one letter in $w$ is a 1, there is a factor in $0w$ of length 1 with more 1s than its prefix, so $s_\sigma[k+1]$ is not prefix normal and $\sigma$ is not a prefix normal word chain generator.\qed
    \end{proof}

\fi

Furthermore, the first flip from $1$ to $0$ cannot happen in the first half of a word. Otherwise, the second half would contain an extra $1$ compared to the first, which also is the prefix of the same length.

\begin{lemma}{}
    \label{lem:firstZeroAtHalf}
    Let $n \in \N$ and let $\sigma$ be a word chain generator of length $n$. If $\sigma$ is a prefix normal word chain generator, then $\sigma(1) \geq \lceil \frac{n+1}{2} \rceil$.
\end{lemma} 
\ifpaper 
\else    
    \begin{proof}
        By contraposition. Let $s_\sigma$ be the sequence of words generated by the word chain generator $\sigma$, and let $w = s_{\sigma}[2]$. Assume that $\sigma(1) < \lceil \frac{n+1}{2} \rceil$. We know that $w = u0v$ for some $u,v \in \{1\}^*$ with 
        $$|u| = \sigma(1) - 1 < \left\lceil \frac{n+1}{2} \right\rceil - 1.$$ 
        And since $|u| + |v| + 1 = |w|$, this means 
        $$|v| = |w| - (|u| + 1) = n - \sigma(1) \geq n - \left\lceil \frac{n+1}{2} \right\rceil + 1.$$
        Now,
        $$ \lceil (n+1) / 2 \rceil - 1 \leq n - \left\lceil \frac{n+1}{2} \right\rceil + 1$$
        gives us $|u| < |v|$, and therefore, $v$ contains one more 1 than the prefix of the same length, $u0\dots$, so $w$ is not prefix normal and $\sigma$ not a prefix normal generator.\qed
    \end{proof}

\fi

For each word length, there is a finite number of prefix normal word chain generators. Often, these are in a relationship with each other: two adjacent numbers in a generator can be swapped, leading to a different prefix normal word chain generator. This word chain will be equal to the first except for one word. A question arises whether all prefix normal word chain generators of a certain length are connected by such swaps. If so, which swaps are allowed, and which lead to the resulting word chain not being prefix normal?

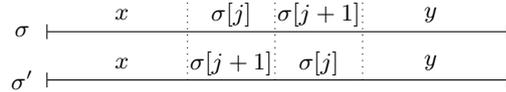
\begin{figure}
    \centering
    \begin{tikzpicture}[scale=0.32,
	blackbox/.style={rectangle, draw=black, thick, minimum height=0.2cm},
	redbox/.style={rectangle, draw=red!75, fill=red!25, thick, minimum height=0.2cm}
	]
    \node at (0,2) {$\sigma$};
    \node at (0,0) {$\sigma'$};
	
    \draw[|-|] (1,2) -- (20,2) node[pos=0.5, sloped] {};	
    \draw[|-|] (1,0) -- (20,0) node[pos=0.5, sloped] {};	

    \node at (4.1, 2.7) {$x$};
    \node at (16.8, 2.7) {$y$};
    \node at (4.1, 0.7) {$x$};
    \node at (16.8, 0.7) {$y$};
    \node at (8.6,0.7) {$\sigma[j+1]$};
    \node at (12.2,2.7) {$\sigma[j+1]$};
    \node at (8.6,2.7) {$\sigma[j]$};
    \node at (12.2,0.7) {$\sigma[j]$};

	\draw[dotted] (6.8, 3.2)--(6.8, -0.2) {};
    \draw[dotted] (10.4, 3.2)--(10.4, -0.2) {};
    \draw[dotted] (14,3.2)--(14,-0.2) {};
    \end{tikzpicture}
    \caption{Prefix normal generator $\sigma$ and \emph{potentially} prefix normal generator $\sigma'$}
    \label{fig:generatorsWithSwap}
\end{figure}

Now, we fix the preliminaries for the following statements. let $\sigma$ be a prefix normal word chain generator of length $n$, and $c_{\sigma}$ the arising prefix normal word chain. Let $j \in [n-1]$, and let $x,y \in \Fact(\sigma)$ such that $\sigma = x \: \sigma[j] \: \sigma[j+1] \: y$. Let $\sigma' = x \: \sigma[j+1] \; \sigma[j] \: y$ (see Figure \ref{fig:generatorsWithSwap}).

\begin{figure}
    \centering
    \begin{tikzpicture}[scale=0.32,
	blackbox/.style={rectangle, draw=black, thick, minimum height=0.6cm},
	redbox/.style={rectangle, draw=red!75, fill=red!25, thick, minimum height=0.6cm}
	]
	\node at (0,6) {$c_{\sigma}[|x|+1]$};
    \node at (0,4) {$c_{\sigma}[|x|+2]$};
    \node at (0,2) {$c_{\sigma'}[|x|+2]$};
    \node at (0,0) {$c_{\sigma}[|x|+3]$};
	
	\draw[|-|] (3,6) -- (20,6) node[pos=0.5, sloped, above] {};
    \draw[|-|] (3,4) -- (20,4) node[pos=0.5, sloped, above] {};	
    \draw[|-|] (3,2) -- (20,2) node[pos=0.5, sloped, above] {};	
    \draw[|-|] (3,0) -- (20,0) node[pos=0.5, sloped, above] {};	
	
    \node at (7.5,6.5) {1};
    \node at (14.5,6.5) {1};
    \node at (7.5,4.5) {0};
    \node at (14.5,4.5) {1};
    \node at (7.5,2.5) {1};
    \node at (14.5,2.5) {0};
    \node at (7.5,0.5) {0};
    \node at (14.5,0.5) {0};

	\draw[dotted] (7,7)--(7,-0.2) {};
    \draw[dotted] (8,7)--(8,-0.2) {};
    \draw[dotted] (14,7)--(14,-0.2) {};
    \draw[dotted] (15,7)--(15,-0.2) {};

    \node at (7.5, -1) {Index $\sigma[j]$};
    \node at (14.5, -1) {Index $\sigma[j+1]$};
	
    \end{tikzpicture}
    \caption{Case $\sigma[j] < \sigma[j+1]$}
    \label{fig:generatorSwapFirstCase}
\end{figure}
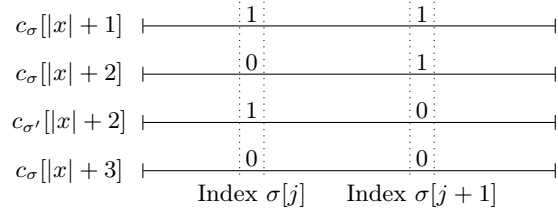

\begin{lemma}{}
    If $\sigma[j] < \sigma[j+1]$, then $\sigma'$ is also a prefix normal word chain generator.
    \label{lem:generatorSwapFirstCase}
\end{lemma}
\ifpaper 
\else    
   \begin{proof}
        Notice that $c_{\sigma'}$ is the new word chain arising from the generator $\sigma'$. The words $c_{\sigma}[|x|+1]$, $c_{\sigma}[|x|+2]$ and $c_{\sigma}[|x|+3]$ are all prefix normal, since they appear in the prefix normal word chain $c_{\sigma}$. We also know that all words in the first word chain $c_{\sigma}$ are the same as in $c_{\sigma'}$ except for one: $c_{\sigma}[i] = c_{\sigma'}[i]$ for $i \in [n], i \neq |x|+1$. The question therefore is whether this word, $c_{\sigma'}[|x|+1]$, is prefix normal (see Figure \ref{fig:generatorSwapFirstCase}).

        We know that $c_{\sigma}[|x|+1]$ contains 1s at its indices $\sigma[j]$ and $\sigma[j+1]$,
        i.e., the $\sigma[j]^{\text{th}}$ letter and the $\sigma[j+1]^{\text{th}}$ letter are both 
        $1$. That is because these two have not yet been flipped to 0, as each number can only appear once in a generator, so they do not occur in $x$. We also know that $c_{\sigma}[|x|+3]$ contains 0s at both indices $\sigma[j]$ and $\sigma[j+1]$. It is only in between these two words that the word chains differ at exactly these two indices: $c_{\sigma}[|x|+2]$ contains a 0 as its $\sigma[j]^{\text{th}}$ letter and a 1 as its $\sigma[j+1]^{\text{th}}$ letter, whereas $c_{\sigma'}[|x|+1]$ contains a 1 at index $\sigma[j]$ and a 0 at index $\sigma[j+1]$.

        Let $u$ be a factor of $c_{\sigma'}[|x|+1]$, which is not a prefix. We will consider several cases based on which of the two letters at index $\sigma[j]$ and $\sigma[j+1]$ are included in $u$. If $u$ contains only the 0, which is the $\sigma[j+1]^{\text{th}}$ letter of $c_{\sigma'}[|x|+1]$, or none of the two letters, then it is also a factor of $c_{\sigma}[|x|+3]$, whose prefixes contain either the same amount or one 1 less than  $c_{\sigma'}[|x|+1]$  which is prefix normal. If $u$ contains only the 1, i.e., the $\sigma[j]^{\text{th}}$ letter of $c_{\sigma'}[|x|+1]$, then it is also a factor of $c_{\sigma}[|x|+1]$, which is prefix normal. There are, of course, prefixes of $c_{\sigma}[|x|+1]$ that contain more 1s than prefixes of the same length of $c_{\sigma'}[|x|+1]$, namely those of length at least $\sigma[j+1]$. These are not relevant here, however, as $u$ must not contain the $\sigma[j+1]^{\text{th}}$ letter in $c_{\sigma'}[|x|+1]$ and therefore is shorter than $\sigma[j+1]$. The last remaining case is if $u$ contains both letters, the 1 and the 0. In that case, $u$ contains the same amount of 1s as the factor starting and ending at the same index in $c_{\sigma}[|x|+2]$. This word is prefix normal and its prefixes contain either the same amount or one 1 less than the prefixes of $c_{\sigma'}[|x|+1]$. So $u$ cannot contain more 1s than the prefix of the same length. Therefore, $c_{\sigma'}[|x|+1]$ is prefix normal.\qed
    \end{proof}

\fi

\begin{figure}
	\centering
	\begin{tikzpicture}[scale=0.32,
		blackbox/.style={rectangle, draw=black, thick, minimum height=0.6cm},
		redbox/.style={rectangle, draw=red!75, fill=red!25, thick, minimum height=0.6cm}
		]
		\node at (0,6) {$c_{\sigma}[|x|+1]$};
		\node at (0,4) {$c_{\sigma}[|x|+2]$};
		\node at (0,2) {$c_{\sigma'}[|x|+2]$};
		\node at (0,0) {$c_{\sigma}[|x|+3]$};
		
		\draw[|-|] (3,6) -- (20,6) node[pos=0.5, sloped, above] {};
		\draw[|-|] (3,4) -- (20,4) node[pos=0.5, sloped, above] {};	
		\draw[|-|] (3,2) -- (20,2) node[pos=0.5, sloped, above] {};	
		\draw[|-|] (3,0) -- (20,0) node[pos=0.5, sloped, above] {};	
		
		\node at (7.5,6.5) {1};
		\node at (14.5,6.5) {1};
		\node at (7.5,4.5) {1};
		\node at (14.5,4.5) {0};
		\node at (7.5,2.5) {0};
		\node at (14.5,2.5) {1};
		\node at (7.5,0.5) {0};
		\node at (14.5,0.5) {0};
		
		\draw[dotted] (7,7)--(7,-0.2) {};
		\draw[dotted] (8,7)--(8,-0.2) {};
		\draw[dotted] (14,7)--(14,-0.2) {};
		\draw[dotted] (15,7)--(15,-0.2) {};
		
		\node at (7.5, -1) {Index $\sigma[j+1]$};
		\node at (14.5, -1) {Index $\sigma[j]$};
		
	\end{tikzpicture} 
	\caption{Case $\sigma[j] > \sigma[j+1]$}
	\label{fig:generatorSwapSecondCase}
\end{figure}
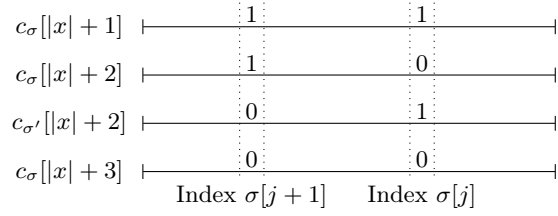

So in this case, $\sigma'$ is always a prefix normal generator. This is in line with our intuition because, starting with the prefix normal word $c_{\sigma}[|x|+2]$, a $1$ is moved closer to the beginning of the word and a $0$ closer to the word's end (see Figure~\ref{fig:generatorSwapFirstCase}). Lemma \ref{lem:generatorSwapFirstCase} also implies that prefix normal words form a bubble language, which has already been shown in \cite{DBLP:conf/cpm/BurcsiFLRS14}.
Now we will move on to the case of $\sigma[j] > \sigma[j+1]$ (Figure~\ref{fig:generatorSwapSecondCase}).

\begin{lemma}{}\label{lem:generatorSwapSecondCase}
    Let $\sigma[j] > \sigma[j+1]$. Let $v$ be a factor of $c_{\sigma}[|x|+1]$ where:\\
    - $v$ includes the  index $\sigma[j]$ and does not include the index $\sigma[j+1],$\\
    - $\sigma[j+1] \leq |v| < \sigma[j]$, and, \\
    - $|v|_1 = |\pref_{|v|}(c_{\sigma}[|x|+1])|_1.$\\
    Then, $\sigma'$ is a prefix normal word chain generator if and only if no such $v$ exists in $c_{\sigma}[|x|+1]$.
\end{lemma}
\ifpaper 
\else    
\begin{proof}
        For the reasons given in the proof of Lemma \ref{lem:generatorSwapFirstCase}, we only have to look at $c_{\sigma'}[|x|+2]$ to determine whether $c_{\sigma'}$ is a prefix normal word chain and, therefore, $\sigma'$ a prefix normal generator. Any factor of $c_{\sigma'}[|x|+2]$ that includes both the 1 at $\sigma[j]$ and the 0 at $\sigma[j+1]$ cannot contain more 1s than the prefix of the same length, since $c_{\sigma}[|x|+2]$ contains a factor of the same length with the same amount of 1s and is prefix normal. Also, any factor of $c_{\sigma'}[|x|+2]$ that does not contain the 1 at index $\sigma[j]$ cannot contain more 1s than the prefix of the same length because it is also a factor of $c_{\sigma}[|x|+3]$, whose prefixes contain either the same amount of 1s as $c_{\sigma'}[|x|+2]$ or one less and which is prefix normal, too. Therefore, any factor of $c_{\sigma'}[|x|+2]$ containing more 1s than the prefix of their length must contain the 1 at index $\sigma[j]$, but not the 0 at $\sigma[j+1]$. 
        
        If the length of such a factor is less than $\sigma[j+1]$, the prefix of the same length will not reach the 0 at that index. In that case, however, both factor and prefix also appear in $c_{\sigma}[|x|+1]$, which is prefix normal. If such a factor is at least $\sigma[j]$ letters long, then the prefix of the same length includes the 1 at index $\sigma[j]$, and both prefix and factor will have one more 1 than the corresponding prefix and factor in $c_{\sigma}[|x|+2]$. And since $c_{\sigma}[|x|+2]$ is prefix normal, adding a 1 to both factor and prefix will also result in a factor with still at most as many 1s as the prefix. This further narrows the factors potentially containing more 1s than the corresponding prefix down to factors of a length of at least $\sigma[j+1]$ and at most $\sigma[j]-1$. 

        We can now use these observations about $c_{\sigma'}[|x|+2]$ and apply them to $c_{\sigma}[|x|+1]$. Let $v$ be a factor of $c_{\sigma}[|x|+1]$ that has all four properties stated in the lemma. This means that there is a factor $v'$ in $c_{\sigma'}[|x|+2]$ at the same position as $v$ with $|v'|_1 = |v|_1$. Because the 1 at index $\sigma[j+1]$ has now been flipped and this 1 is included in the prefix of length $|v'|$, we have $|\pref_{|v|}(c_{\sigma}[|x|+1])|_1 - 1 = |\pref_{|v'|}(c_{\sigma'}[|x+2|])|_1$. This means that $v'$ contains exactly one more 1 than the prefix of the same length. So $c_{\sigma'}[|x|+1]$ is not prefix normal, and the swap results in $c_{\sigma'}$ not being a prefix normal generator. 
        
        To show that, if $\sigma'$ is a prefix normal word chain generator, no such $v$ can exist, it suffices to show the following: If a $v$ fulfilling the properties exists, then $\sigma'$ is not a prefix normal word chain generator. If such a $v$ exists in $c_{\sigma}[|x|+1]$, then we know the following about $c_{\sigma'}[|x|+2] = uvw$ for some $u, w \in \Sigma^*$: $|u| > 0$, because $v$ excludes the letter at index $\sigma[j+1]$ but includes the letter at index $\sigma[j]$, and $\sigma[j] > \sigma[j+1]$. We also know $|v|_1 > |\pref_{|v|}(c_{\sigma'}(|x|+2)|_1$, because this prefix now contains a 0 at index $\sigma[j+1]$ which $v$ does not contain, and both are otherwise unchanged when compared to the factors at the same indices in $c_{\sigma}(|x|+1)$, where $v$ contained the same amount of 1s as the prefix of the same length. Together, these two observations give us that $v$ is a factor that is not a prefix and contains more 1s than the prefix of the same length. So $c_{\sigma'}(|x|+2)$ is not a prefix normal word and hence $\sigma'$ is not a prefix normal word chain generator.
        \qed
\end{proof}
\fi

\begin{remark}
	In the light of Lemma~\ref{lem:generatorSwapSecondCase}, the question arises how many possible factors one needs to test for deciding whether $\sigma'$ is a prefix normal word chain generator. 
	The worst case on the number of possible factors $v$ that need to be tested on the conditions of Lemma~\ref{lem:generatorSwapSecondCase} happens for $\sigma[j+1] = 2$ and $\sigma[j]=\lfloor\frac{n-2}{2}\rfloor$ ($\sigma[j+1]=1$ is not relevant by Lemma~\ref{lem:1ComesLastInGen}).
	By the first condition ($v$ includes index $\sigma[j]=\lfloor\frac{n-2}{2}\rfloor$ and does not include index $\sigma[j+1]=2$) we know that there exist $\lfloor\frac{n-2}{2}\rfloor$ such factors of length $\lfloor\frac{n-2}{2}\rfloor$, $\lfloor\frac{n-2}{2}\rfloor -1$ such factors of length $\lfloor\frac{n-2}{2}\rfloor -1$ etc. Thus, the number of factors is upper bounded by the Gaussian sum from $1$ to $\lfloor\frac{n-2}{2}\rfloor$, i.e., $\frac{\frac{n-2}{2}(\frac{n-2}{2} +1)}{2} = \frac{n^2-2n}{8}$.
	
	Note, that the largest possible range in the condition $\sigma[j+1] \leq |v| < \sigma[j]$ of Lemma~\ref{lem:generatorSwapSecondCase} is given by $\sigma[j+1] = 2$ and $\sigma[j]=n$ does not automatically lead to the worst case. By the first condition ($v$ includes index $\sigma[j]=n$ and does not include index $\sigma[j+1]=2$) we know that every suffix of $c_\sigma[|x|+1]$ from the $3^\text{rd}$ position on needs to be inspected (the suffix of length $1$ is  not relevant since we ensured that $\sigma(n) =1$ holds for prefix normal word chain generator), i.e., less factors than in the previously described case.
	
	Further, one can apply the results from Section~\ref{sec:minFactors} in order to reduce the number of inspected factors even more (e.g., only inspecting those that start and end with the letter $1$). 
 \end{remark}

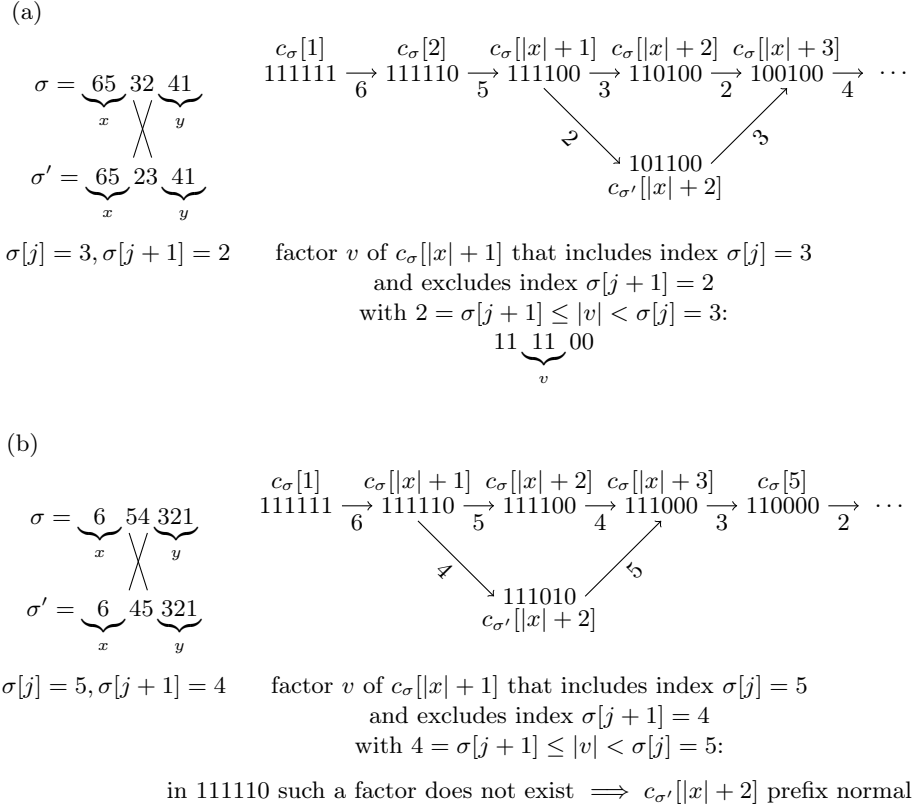
\begin{figure}[p]
	\centering
	\begin{tikzpicture}[scale=0.4,
		blackbox/.style={rectangle, draw=black, thick, minimum height=0.6cm},
		redbox/.style={rectangle, draw=red!75, fill=red!25, thick, minimum height=0.6cm}
		]
		
		\node at (-3,9) {(a)};
		
		\node at (0,6) {$\sigma = \underbrace{65}_{x}32\underbrace{41}_{y}$};
		\node at (0,3) {$\sigma' = \underbrace{65}_{x}23\underbrace{41}_{y}$};
		\node at (0,1) {$\sigma[j] = 3, \sigma[j+1] = 2$};
		
		\draw[-] (0.5,6.1) -- (1.1,4.2) node[pos=0.5, sloped, above] {};
		\draw[-] (1.1,6.1) -- (0.5,4.2) node[pos=0.5, sloped, above] {};
		
		\draw[->] (7.5,7) -- (8.5,7) node[pos=0.5, sloped, below] {6};	
		\draw[->] (11.5,7) -- (12.5,7) node[pos=0.5, sloped, below] {5};	
		\draw[->] (15.5,7) -- (16.5,7) node[pos=0.5, sloped, below] {3};	
		\draw[->] (19.5,7) -- (20.5,7) node[pos=0.5, sloped, below] {2};	
		\draw[->] (23.5,7) -- (24.5,7) node[pos=0.5, sloped, below] {4};	
		
		\draw[->] (14,6.5) -- (16.5,4) node[pos=0.5, sloped, below] {2};	
		\draw[->] (19.5,4) -- (22,6.5) node[pos=0.5, sloped, below] {3};	
		
		\node at (6,7.8) {$c_\sigma[1]$};
		\node at (6,7) {111111};
		
		\node at (10,7.8) {$c_\sigma[2]$};
		\node at (10,7) {111110};
		
		\node at (14,7.8) {$c_\sigma[|x|+1]$};
		\node at (14,7) {111100};
		
		\node at (18,7.8) {$c_\sigma[|x|+2]$};
		\node at (18,7) {110100};
		
		\node at (22,7.8) {$c_\sigma[|x|+3]$};
		\node at (22,7) {100100};
		
		\node at (25.5,7) {$\ldots$};
		
		\node at (18,3.2) {$c_{\sigma'}[|x|+2]$};
		\node at (18,4) {101100};		
		
		\node at (14,1) {factor $v$ of $c_\sigma[|x|+1]$ that includes index $\sigma[j] = 3$};
		\node at (14,0) {and excludes index $\sigma[j+1] = 2$};
		\node at (14,-1) {with $2= \sigma[j+1] \leq |v| < \sigma[j] =3$:};
		\node at (14,-2.5) {$11\underbrace{11}_{v}00$};
	\end{tikzpicture} 
	\\~\\
	\begin{tikzpicture}[scale=0.4,
		blackbox/.style={rectangle, draw=black, thick, minimum height=0.6cm},
		redbox/.style={rectangle, draw=red!75, fill=red!25, thick, minimum height=0.6cm}
		]
		\node at (-3,9) {(b)};
		
		\node at (0,6) {$\sigma = \underbrace{6}_{x}54\underbrace{321}_{y}$};
		\node at (0,3) {$\sigma' = \underbrace{6}_{x}45\underbrace{321}_{y}$};
		\node at (0,1) {$\sigma[j] = 5, \sigma[j+1] = 4$};
		
		\draw[-] (0.5,6.1) -- (1.1,4.2) node[pos=0.5, sloped, above] {};
		\draw[-] (1.1,6.1) -- (0.5,4.2) node[pos=0.5, sloped, above] {};
		
		\draw[->] (7.5,7) -- (8.5,7) node[pos=0.5, sloped, below] {6};	
		\draw[->] (11.5,7) -- (12.5,7) node[pos=0.5, sloped, below] {5};	
		\draw[->] (15.5,7) -- (16.5,7) node[pos=0.5, sloped, below] {4};	
		\draw[->] (19.5,7) -- (20.5,7) node[pos=0.5, sloped, below] {3};	
		\draw[->] (23.5,7) -- (24.5,7) node[pos=0.5, sloped, below] {2};	
		
		\draw[->] (10,6.5) -- (12.5,4) node[pos=0.5, sloped, below] {4};	
		\draw[->] (15.5,4) -- (18,6.5) node[pos=0.5, sloped, below] {5};	
		
		\node at (6,7.8) {$c_\sigma[1]$};
		\node at (6,7) {111111};
		
		\node at (10,7.8) {$c_\sigma[|x|+1]$};
		\node at (10,7) {111110};
		
		\node at (14,7.8) {$c_\sigma[|x|+2]$};
		\node at (14,7) {111100};
		
		\node at (18,7.8) {$c_\sigma[|x|+3]$};
		\node at (18,7) {111000};
		
		\node at (22,7.8) {$c_\sigma[5]$};
		\node at (22,7) {110000};
		
		\node at (25.5,7) {$\ldots$};
		
		\node at (14,3.2) {$c_{\sigma'}[|x|+2]$};
		\node at (14,4) {111010};		
		
		\node at (14,1) {factor $v$ of $c_\sigma[|x|+1]$ that includes index $\sigma[j] = 5$};
		\node at (14,0) {and excludes index $\sigma[j+1] = 4$};
		\node at (14,-1) {with $4= \sigma[j+1] \leq |v| < \sigma[j] =5$:};
		\node at (14,-2.5) {in $111110$ such a factor does not exist $\implies c_{\sigma'}[|x|+2]$ prefix normal};
	\end{tikzpicture}
	\caption{Examples for the factor $v$ from Lemma~\ref{lem:generatorSwapSecondCase}: (a) $\sigma = 653241 \to \sigma' = 652341$ is not a pnwcg; (b) $\sigma = 654321 \to \sigma' = 645321$ is a pnwcg.}
	\label{fig:examplesforv}
\end{figure}

An extended example on the existence of a factor $v$ as in Lemma~\ref{lem:generatorSwapSecondCase} can be found in Figure~\ref{fig:examplesforv}. We will now combine these last two lemmas into one general observation about the prefix normality of generators after a swap.

\begin{theorem}{}
		 A generator $\sigma'$ is not a prefix normal generator if and only if \\
			- $\sigma[j] > \sigma[j+1]$, and\\
			- a factor $v \in  \Fact(c_{\sigma}[|x|+1])$ exists such that \\
			\hphantom{factor} - $v$ includes index $\sigma[j]$ and does not include index $\sigma[j+1],$\\
			\hphantom{factor} - $\sigma[j+1] \leq |v| < \sigma[j]$ and\\
			\hphantom{factor} - $|v|_1 \geq |\pref_{|v|}(c_{\sigma}[|x|+1])|_1.$
\end{theorem}
\ifpaper 
\else    
    \begin{proof}
    See Lemma \ref{lem:generatorSwapFirstCase} and Lemma \ref{lem:generatorSwapSecondCase}.\qed
        
        
    \end{proof}

\fi

This means that we do not have to look at $c_{\sigma'}[|x|+2]$ or indeed any part of $c_{\sigma'}$ itself to determine whether it is prefix normal. It is sufficient to check a prefix normal word chain generator by looking at pairs of entries like this: If the first of the two numbers is smaller than the second, swapping them will lead to another prefix normal generator. If the first number is larger than the second, we can test the last word in the corresponding word chain --- where all flips up to the two that the number pair specifies, have already happened --- for factors with the specified characteristics.


Instead of swapping elements of prefix normal generators, we can also increase their length. Prefix normal word chains and generators differ from the way prefix normal words have been studied in that the words we consider are all of the same length \cite{DBLP:conf/fun/BurcsiFLRS14,DBLP:journals/tcs/BurcsiFLRS17}. Now, we combine this view with the classical approach of appending letters to words. First, we will look at complete graphs and ask the question: What is the difference between a graph showing all prefix normal words of length $n$ and its successor showing all prefix normal words of length $n+1$? First we give some easy results we need to proceed.

\begin{remark}{(cf. \cite{DBLP:conf/dlt/FiciL11})}
    If a word $w \in \Sigma^*$ is prefix normal, so are $1w$ and $w0$.
    \label{rem:Prepend1Append0}
\end{remark}

\begin{lemma}{}
    The number of prefix normal words of length $n$ is equal to the number of prefix normal words of length $n+1$ ending in 0.
\end{lemma}
\ifpaper 
\else    
  \begin{proof}
        Let $n, m \in \N$, $\text{PNW}_n = \{w \: | \: w \text{ prefix normal and } |w|=n\}$, and let $\text{PNW0}_m = \{w \: | \: w \text{ prefix normal, } |w|=m \text{ and } w[m] = 0\}$. We want to show $|\text{PNW0}_{n+1}| = |\text{PNW}_n|$.
        
        
        Firstly, note that if a 0 is appended to a prefix normal word, the resulting word is also prefix normal (see Remark \ref{rem:Prepend1Append0}). This means that by appending a 0 to all prefix normal words of length $n$, we get a set of prefix normal words ending in 0. So we have $|\text{PNW0}_{n+1}| \geq |\text{PNW}_n|$.
        
        Now, assume that $|\text{PNW0}_{n+1}| > |\text{PNW}_n|$. This means there exists a word $w \in \Sigma^*$ with $|w| = n+1, w[n+1] = 0$, which is prefix normal, but $w[1:n]$ is not prefix normal. So $w[1:n]$ contains a factor which includes more 1s than the prefix of the same length. Since both prefix and factor also occur in $w$, and $w$ is prefix normal, this cannot be the case. So $w[1:n]$ is prefix normal, and therefore, $|\text{PNW0}_{n+1}| \leq |\text{PNW}_n|$.

        Thus we have $|\text{PNW0}_{n+1}| = |\text{PNW}_n|$.\qed
    \end{proof}

\fi

So we know more than half of the prefix normal words of length $n + 1$ by appending a 0 to all words of length $n$. But how do we find the remaining words? One way to look at it is by asking: To which words of the graph of prefix normal words of length $n$ can we append a 1 with the resulting word being prefix normal? And what are the characteristics of words where this is not possible? In \cite{DBLP:conf/fun/BurcsiFLRS14,BURCSI202086}, a prefix normal word $w$ is called \emph{extension-critical} if $w1$ is not prefix normal.

\begin{lemma}{}
    If $1$ is appended to a prefix normal palindrome containing at least one $0$, the resulting word is not prefix normal.
    \label{lem:PalindromeExtensionCritical}
\end{lemma}
\ifpaper 
\else    
    \begin{proof}
        Let $w \in \Sigma^*$, $w$ prefix normal and a palindrome, $|w|_0 > 0$. Let $x \in \Pref(w)$ such that $|x|_0 = 0$ and $x0 \in \Pref(w)$. As $w$ is a palindrome and $x$ contains only 1s, we know $0x \in \Suff(w)$. This means that $x0 \in \Pref(w1)$ and $x1 \in \Suff(w1)$, and since $|x0|_1 = |x|_1 < |x1|_1$, $w1$ is not prefix normal.\qed
    \end{proof}

\fi

Going one step further, we can combine prefix normal word chain generators with the approach of appending letters like this: Take a prefix normal generator of length $n$, and insert the number $n+1$ at different points. This means appending one letter to each word in the word chain. In which cases will the new generator (not) be prefix normal? It follows immediately from Lemma \ref{lem:1ComesLastInGen} that $n+1$ cannot be inserted at the last position of the new generator. Let $\sigma$ be a prefix normal generator of length $n$ and $c_{\sigma}$ its corresponding prefix normal word chain. Let $w = c_{\sigma}[i]$ for some $i \in \N$. Let $\tau$ be the word chain generator arising from the insertion of $n+1$ into $\sigma$ such that $w$ is the last word in $c_{\sigma}$ to which $1$ is appended, while $0$ is appended to its child $c_{\sigma}[i+1]$ and all other words that follow in the word chain. So
$$\tau = (\sigma[1], \sigma[2], ..., \sigma[i-1], n+1, \sigma[i], ..., \sigma[n])\text{, and}$$
$$c_{\tau} = (c_{\sigma}[1]1, c_{\sigma}[2]1, ..., w1, c_{\sigma}[i+1]0, ..., c_{\sigma}[n+1]0, 0^{n+1}).$$

Now, if $w1$ is not prefix normal, neither is $c_{\tau}$. But even if $w1$ \emph{is} prefix normal, there are cases where $c_{\tau}$ is not. What we want to show is that, if $w1$ is prefix normal, it is always part of a prefix normal word chain that was constructed in the same way as $c_{\tau}$. We know from Remark \ref{rem:Prepend1Append0} that $c_{\sigma}[i+1]0$ and all words following it in $c_{\tau}$ are prefix normal. So we only have to look at the first $i$ words in $c_{\tau}$. We do this by showing that every word like $w1$ and the words before it in $c_{\tau}$ have at least one prefix normal parent. If we look at all prefix normal words of length $n \in \N$ as a graph, this means the following: We append 1 to every word in the graph that contains less than $k \in \N$ 0s, and we append 0 to every word in the graph containing at least $k$ 0s. Then, we take a prefix normal word containing exactly $k-1$ occurrences of 0 from the resulting graph. This word will be connected to the first and last word in the graph (the two words containing only 1s and 0s) via a word chain of only prefix normal words.

\begin{proposition}{}
    Let $w$ be a prefix normal word. If $w1$ is also prefix normal, then $w$ has at least one prefix normal parent $v$ such that $v1$ is also prefix normal.

    \label{lem:NoIsolatedApp1}
\end{proposition}
\ifpaper 
\else    
    \begin{proof}
        We have to show that such a $v$ exists. Let $v$ be the parent of $w$ that has a longer prefix containing only 1s, i.e., for $w = 1^k0x, k \in \N, x \in \Sigma^*$, let $v = 1^{k+1}x$. Assume that $w1$ is prefix normal, and let $v$ be a prefix normal parent of $w$. Any factor in $v1$ that might lead to it not being prefix normal would have to include both the 1 at the end as well as the 1 at position $(k+1)$. If the last 1 is not included, the factor and the corresponding prefix already occur in $v$, which is prefix normal. If the one at position $k+1$ is not included, the factor already occurs in $w1$, and the corresponding prefix either has the same amount of 1s as the corresponding prefix in $w1$ or one more. However, the 1 at the $(k+1)^{\text{th}}$ position is part of the prefix containing only 1s, and Lemma \ref{lem:FacInPrefixOfOnes} tells us that a factor responsible for a word not being prefix normal cannot begin inside the prefix of that word which contains only 1s. Therefore, there is no factor containing more 1s than the prefix of the same length in $v1$, so $v1$ is prefix normal.\qed
    \end{proof}

\fi

Further, we know that for a word chain generator that is not prefix normal, a right shift of the entry $n$ always results again in a word chain generator that is not prefix normal.

\begin{proposition}\label{prop:ShiftOfHighestValueStaysNonPrefixnormal}
	If $\sigma'$ is a prefix normal word chain generator of length $n \in \N$ with $\sigma'[i] = n$ for some $i \in [2,n-1]$ then $\sigma$ given by $\sigma = \sigma[1..i-2] \sigma[i] \sigma[i-1]  \sigma[i+1..n]$ is a prefix normal word chain generator.
\end{proposition}
\ifpaper 
\else    
\begin{proof}
	We are showing the contraposition of the claim: If $\sigma$ is a not prefix normal word chain generator of length $n$ and $\sigma[i] = n$ for some $i \in [n-1]$ then every $\sigma'$ with $\sigma' = \sigma[1..i-1]\sigma[i+1] \sigma[i] \sigma[i+2..n]$ is not a prefix normal word chain generator.
	
	First, we observe that there are three options: either the first prefix normal word within the word chain $c_\sigma$ is at a position $j<i, j=i$ or $j>i$.
	
	If the first non prefix normal word is generated at a position $j$ in $c_\sigma$ that is smaller than $i$, we immediately obtain that every $\sigma' = \sigma[1..i-1]\sigma[i+1] \sigma[i] \sigma[i+2..n]$ is not a prefix normal word chain generator witnessed by the non prefix normal word at position $c_\sigma[j] = c_{\sigma'}[j]$. Thus, $c_{\sigma'}$ is not a prefix normal word chain generator.
	
	If the first non prefix normal word in the word chain $\sigma$ is at position $i$. Thus, $w = c_\sigma[i-1]$ is a prefix normal word. Since $w$ is prefix normal, we can also conclude that $w[1..n-1]0$ is prefix normal. Since $w[1..n-1]0 = c_\sigma[i]$, we obtain a contradiction to the fact that the first non prefix normal word in the word chain $\sigma$ is at position $i$, which means that this case cannot happen.
	
	If now $j > i$, we know that $c_\sigma[i] = x0$ for some $x \in \Sigma^*$ is prefix normal. Further, $c_\sigma[j] = x' 0 y' 0$ is the first prefix normal word within the word chain (preceded by $c_\sigma[j-1] = x'1 y' 0$) for some $x',y' \in \Sigma^*$. This means, that there exists a factor of $x' 0 y' 0$ that contains more $1$s than the prefix $x'0$. Since this factor does not contain $c_\sigma[n]$ we can conclude that the position of $n$ within the word chain generator does not influence $c_\sigma$ of not being prefix normal chain generator, i.e., $c_{\sigma'}$ is not a prefix normal word chain generator.\qed
\end{proof}
\fi

Using the previous results, we give the idea of constructing the prefix normal generators of length $n+1$ from those of length $n$ the following way: We take a prefix normal generator of length $n$ and insert $n+1$ as the new first element of the sequence. This means creating the new prefix normal word $1^{n+1}$ as the new first element of the word chain and appending a 0 to all other words in it. We know that the resulting generator is prefix normal. Now, we move the number $n+1$ one position to the right, i.e., we swap it with the element succeeding it in the sequence. If the resulting generator is prefix normal, we add it to the list of prefix normal generators of length $n+1$ and repeat this step. If not, we stop (by contraposition of \Cref{prop:ShiftOfHighestValueStaysNonPrefixnormal}) and continue with the next generator of length $n$. From \Cref{lem:NoIsolatedApp1}, we know that we are not missing prefix normal words of length $n+1$ ending in 1. That is because every prefix normal word to which a 1 can be appended is connected to a parent to which a 1 can also be appended. So every such prefix normal word will be found because it is part of a word chain where all previous words also remain prefix normal if a 1 is appended. This algorithm can also include a check whether the word at the current position of the generator of length $n$ is a palindrome, in which case no 1 can be appended because of \Cref{lem:PalindromeExtensionCritical}.

To find all prefix normal words of a certain length, this algorithm has to be applied recursively, with the prefix normal generator of length 1 being given as (1). It has to test almost every prefix normal word of length $n$ on whether it is extension-critical to find those of length $n+1$. Its complexity therefore is too high for it to be a practical approach. However, as stated similarly in the conclusion of \cite{DBLP:conf/cpm/BurcsiFLRS14}, the number of binary words grows much more rapidly than that of prefix normal words. So if all prefix normal words of length $n$ are given, the algorithm performs significantly better than generating all binary words of length $n+1$ and then filtering on the condition given in the definition of prefix normality.

Since permutations and the symmetric group are well-studied objects, we investigated two basic notions of permutations. First, the closure of composition within the symmetric group and second, the parity of a permutation. 

\begin{remark}
	Note that the composition of two prefix normal generators is never a prefix normal generator again. Let $\sigma,\sigma'$ be prefix normal generators over $[n]$ for $n \in \N$. By \Cref{lem:1ComesLastInGen} we have that $\sigma(n) = \sigma'(n) =1$. Thus, in the composition $\sigma \circ \sigma'$ we get $(\sigma \circ \sigma')(n)\sigma(\sigma'(n)) = \sigma(1) \neq 1$, i.e., the composition cannot be a prefix normal generator.
\end{remark}

\begin{remark}
	There exist prefix normal generators with even and odd parity witnessed by $(3241), (3421), (4321), (4231)$ that have $4,5,6,5$ inversions, so even, odd, even, odd parities, resp.
\end{remark}

To conclude this section, we presented a new perspective on prefix normal words of a fixed length $n$ giving a characterisation of prefix normal word chains. For now a connection to existing results in theory about permutations and the symmetric group cannot be made and stay open for further research.


\section{Conclusion}\label{conc}
Our work on prefix normal words is based on \cite{DBLP:conf/dlt/FiciL11}. In  \Cref{sec:minFactors}, we presented characteristics of factors that contain more 1s than the prefix of the same length, hence exploiting the word that is responsible for them being not prefix normal. 

Our work in Section 4 focused on relationships between words of equal length. We introduced word chains and word chain generators that describe multiple words of the same length. While these concepts enabled us to think about prefix normal words in a more abstract way, the problem of enumerating prefix normal words is open to further investigation. Research into both the relationship between generators of the same length as well as inserting one number into a generator is necessary. We have shown that testing whether a swap will lead to a new prefix normal generator can be done by testing only specific factors of a word. Is there a path from any prefix normal generator to all others, just by swapping neighbouring numbers? If so, all prefix normal words of a length could be found by checking possible swap positions. It is not yet clear whether these positions follow a pattern or are further related to the symmetric group.
In \cite{DBLP:conf/dlt/EikmeierFKN21}, \emph{weighted} prefix normal words were introduced, extending the concept of prefix normality to arbitrary finite alphabets. Word chain generators could become \emph{word chain generating matrices} with $n \in \N$ rows to represent word chains of words over alphabets containing $n+1$ letters.


\bibliographystyle{plainurl}
\bibliography{refs}



\end{document}